\documentclass{amsart}

\usepackage{latexsym}
\usepackage{amsmath}
\usepackage{amssymb}
\begin{document}


\newcommand{\az}{\alpha_{Z}}
\newcommand{\bz}{\beta_{Z}}
\newcommand{\ax}{\alpha_{\X}}
\newcommand{\bx}{\beta_{\X}}
\newcommand{\popo}{\mathbb{P}^1 \times \mathbb{P}^1}
\newcommand{\pnk}{\mathbb{P}^{n_1}\times \cdots \times \mathbb{P}^{n_k}}
\newcommand{\Iz}{I_{Z}}
\newcommand{\Ix}{I_{\X}}
\newcommand{\C}{\mathcal{C}}
\newcommand{\Y}{\mathbb{Y}}
\newcommand{\Z}{\mathbb{Z}}
\newcommand{\Zr}{\mathbb{Z}_{red}}
\newcommand{\N}{\mathbb{N}}
\newcommand{\pr}{\mathbb{P}}
\newcommand{\X}{\mathbb{X}}
\newcommand{\supp}{\operatorname{Supp}}
\newcommand{\ol}{\overline{L}}
\newcommand{\Ssz}{\mathcal S_{Z}}
\newcommand{\Ss}{\mathcal S}
\newcommand{\dtz}{\Delta H_{Z}}
\newcommand{\dtc}{\Delta^{C} H_{Z}}
\newcommand{\dtcc}{\Delta^{C} H_{Z_{ij}}}
\newcommand{\dt}{\Delta}
\newcommand{\B}{\mathcal{B}}
\newcommand{\ay}{\alpha_{Y}}
\newcommand{\by}{\beta_{Y}}
\newcommand{\Q}{\mathcal{Q}}

\newtheorem{theorem}{Theorem}[section]
\newtheorem{corollary}[theorem]{Corollary}
\newtheorem{proposition}[theorem]{Proposition}
\newtheorem{lemma}[theorem]{Lemma}
\newtheorem{alg}{Algorithm}
\newtheorem{question}{Question}
\newtheorem{problem}{Problem}
\newtheorem{conjecture}[theorem]{Conjecture}

{\theoremstyle{definition}
\newtheorem{remark}[theorem]{Remark}
\newtheorem{convention}[theorem]{Convention}
\newtheorem{definition}[theorem]{Definition}
\newtheorem{example}[theorem]{Example}
}


\title{Fat Points in $\popo$ and their Hilbert functions}
\thanks{Revised Version: August 6, 2002}
\author{Elena Guardo}
\address{Dipartimento di Matematica e Informatica\\
Viale A. Doria, 6 - 95100 - Catania, Italy}
\email{guardo@dmi.unict.it}
\author{Adam Van Tuyl}
\address{Department of Mathematics \\
Lakehead University \\
Thunder Bay, ON P7B 5E1, Canada}
\email{avantuyl@sleet.lakeheadu.ca}
\keywords{Hilbert function, points, fat points, Cohen-Macaulay
multi-projective space}
\subjclass{13D40,13D02,13H10,14A15}

\begin{abstract}
We study the Hilbert functions of fat points in $\popo$.
If $Z \subseteq \popo$ is an arbitrary fat point scheme, then
it can be shown that for every $i$ and $j$ the values of the Hilbert function
$H_{Z}(l,j)$ and $H_{Z}(i,l)$ eventually become constant for
$l \gg 0$.  We show how to determine these eventual values
by using only the multiplicities of the points, and the
relative positions of the points in $\popo$.   This enables
us to compute all but a finite number values of $H_{Z}$
without using the coordinates of points.
We also  characterize the ACM fat points schemes
using our description of the eventual behaviour.  In fact,
in the case  that $Z \subseteq \popo$ is ACM, then
the entire Hilbert function and its minimal free resolution
depend solely on knowing the eventual values of the Hilbert function.
\end{abstract}

\maketitle


\section*{Introduction}
The Hilbert function of a fat point scheme in $\pr^n$ is the
basis for many questions about fat points schemes.  Although
some facts have been established (see the survey of Harbourne
\cite{H} for the case of $n=2$), we do not have a complete
understanding of the Hilbert functions of fat point schemes.

In this paper we investigate the Hilbert functions
of fat point schemes in a different space, specifically,
in $\popo$.  Interest in the Hilbert functions
of fat point schemes in $\pnk$ with $k \geq 2$
is motivated, in part, by the work of Catalisano, {\it et al.}
\cite{CGG} which exhibited a connection
between a specific value of the Hilbert function
of a special fat point scheme in $\pnk$ and a classical
problem of computing the dimension
of certain secant varieties to the Segre variety.

The Hilbert functions of sets of points in $\popo$ appear to be
first studied by Giuffrida, {\it et al.} \cite{GuMaRa1}.  Some
of the results of ~\cite{GuMaRa1} were extended and generalized
to sets of points in $\pnk$ by the second author \cite{VT1,VT2}.
Unlike the case of sets of simple points in $\pr^n$, the
problem of characterizing the Hilbert functions of sets of
reduced points in $\pnk$, even in the case of $\popo$, remains
open.  Arithmetically Cohen-Macaulay fat point schemes in $\popo$ were
studied by the first author ~\cite{Gu} (which
was based upon ~\cite{Gu2}).  Catalisano, {\it et al.} ~\cite{CGG}
give some results about fat point schemes in $\pnk$.  However, like the
case of fat point schemes in $\pr^n$, we do not have a complete understanding
of the Hilbert functions of fat point schemes in $\pnk$.

In this paper we are specifically interested in studying
the eventual behaviour of the Hilbert function
of a fat point scheme $Z \subseteq \popo$.
If $Z$ is an arbitrary fat point scheme and if $H_{Z}$
denotes its Hilbert function, then it is not
difficult to show that for any $i$ or $j$, the values
$H_{Z}(l,j)$ and $H_{Z}(i,l)$ become constant for $l \gg 0$.
Our first main result (Theorem \ref{eventualbehaviour})
is to calculate these eventual values
by using numerical information about $Z$.
In particular, we show that these values
can be calculated directly from the multiplicities of the points, and
from the relative positions of the points in the support, that is,
if $P,P'$ are in the support, we only need to know
if $\pi_i(P) = \pi_i(P')$ for $i = 1,2$ where $\pi_i$
is the $i$-th projection map.  The actual coordinates
of the points are therefore not needed to compute
all but a finite number of values of $H_{Z}$.

We then show that the eventual behaviour of $H_{Z}$
gives us further information about
the scheme $Z$.
In particular, we show (cf. Theorem \ref{equivalent})
that the eventual values of $H_{Z}$ can be
used to determine if $Z$ is arithmetically Cohen-Macaulay (ACM).
In fact, a specific type of eventual behaviour
characterizes the ACM fat point schemes of $\popo$.
We relate our characterization with the results
of \cite{GuMaRa1} and \cite{Gu}.
Furthermore,
in the case that $Z$ is ACM, the eventual values
of  $H_{Z}$ can be used to completely determine the entire
Hilbert function, {\it and} the minimal free resolution, of $Z$.

This paper has five parts.  In the first section we recall the
relevant facts about bigraded rings and fat point schemes.  We
also give some elementary properties for the Hilbert function of
a fat point scheme in $\popo$.  In the second section we compute
the Hilbert function of a fat point scheme in $\popo$ whose
support lies on either a $(0,1)$-line or a $(1,0)$-line. In the
third section we introduce two tuples $\az$ and $\bz$
that contain information about the multiplicities and relative
position of the points, and show
how to compute all but a finite number of values of the Hilbert
function from $\az$ and $\bz$.  In the fourth section we show
how to use $\az$ and $\bz$ to determine if $Z$ is ACM.  In the
final section, we look at some ACM fat point schemes with some
extra conditions on their multiplicities.

Many of these results had their genesis in examples.  Instrumental
in computing these examples was the computer program {\tt CoCoA} ~\cite{C}.
We would like to thank A. Ragusa for his useful comments and
suggestions.  We would also like to thank the referee
for their helpful comments and suggestions, and especially
for suggesting a shorter proof for
Theorem \ref{fatpointsonaline}.


\section{Preliminaries}

In this section we recall the necessary definitions and facts about
bigraded rings and fat point schemes.

Let $\N := \{0,1,2,\ldots\}$.
It will be useful to consider in $\Bbb {Z \times
Z}$ and in $\Bbb {N\times \dots \times N}$ the partial
 ordering induced by the usual one in $\Bbb Z$ and
in $\Bbb N$ respectively. We will denote it by \lq\lq
$\leq$\rq\rq.  Thus, if  $(i_1,i_2),(j_1,j_2) \in \N^2$, then we
write $(i_1,i_2) \leq (j_1,j_2)$ if $i_k \leq j_k$ for $k =1,2$.

We let ${\bf k}$ denote an algebraically closed field.
Let $R ={\bf k}[x_0,x_1,y_0,y_1]$ where
$\deg x_i = (1,0)$ and $\deg y_i = (0,1)$.  Then the ring $R$ is
$\N^2$-{\it graded}, or simply, {\it bigraded}, that is,
\[
R = \bigoplus_{(i,j) \in \N^2} R_{i,j}
\hspace{.5cm}\text{and}\hspace{.5cm} R_{i_1,i_2}R_{j_1,j_2}
\subseteq R_{i_1+j_1,i_2+j_2}
\]
were each $R_{i,j}$ consists of all the {\it bihomogeneous
elements} of degree $(i,j)$.

For each $(i,j) \in \N^2$, the set $R_{i,j}$ is a finite
dimensional vector space over ${\bf k}$.  A basis for $R_{i,j}$
is the set of monomials
$
\{x_0^{a_0}x_1^{a_1}y_0^{b_0}y_1^{b_1} \in R ~|~
(a_0+a_1,b_0+b_1) = (i,j)\}.
$
It follows that  $\dim_{\bf k} R_{i,j} = (i+1)(j+1)$ for all
$(i,j) \in \N^2.$

Suppose that $I = (F_1,\ldots,F_r) \subseteq R$ is an ideal such
that the $F_i$'s are bihomogeneous elements.  Then $I$ is called
a {\it bihomogeneous ideal}.  If $I \subseteq R$ is any ideal,
then we define $I_{i,j} := R_{i,j} \cap I$.  The set $I_{i,j}$
is  a subvector space of $R_{i,j}$.  If $I$ is a bihomogeneous
ideal, then $I = \bigoplus_{(i,j)} I_{i,j}$.

If $I$ is a bihomogeneous ideal of $S$, then the quotient ring $S
= R/I$ is also bigraded, i.e., $S = \bigoplus_{(i,j)} S_{i,j}$
where $S_{i,j} := R_{i,j}/I_{i,j}$ for all $(i,j) \in \N^2$.  The
numerical function $H_{S}:\N^2 \rightarrow \N$ defined by
\[
(i,j) \longmapsto \dim_{\bf k} S_{i,j} = \dim_{\bf k} R_{i,j} -
\dim_{\bf k} I_{i,j}
\]
is the {\it Hilbert function of $S = R/I$}.  We sometimes write
the values of the Hilbert function $H_S$ as an infinite matrix
$(M_{i,j})$ where $M_{i,j}:= H_S(i,j)$.  For example, if $I =
(0)$, then $H_{R/I}(i,j) = (i+1)(j+1)$, and so we write
\[
H_{R/I} =
\bmatrix
1 &2&3&4 &\cdots \\
2&4&6&8&\cdots \\
3&6&9&12&\cdots \\
4&8&12&16&\cdots \\
\vdots&\vdots&\vdots&\vdots&\ddots \\
\endbmatrix.
\]
Note that we begin the indexing of the rows and columns at $0$
rather than $1$.

\begin{remark}
In \cite{GuMaRa1} the Hilbert function was referred to as the Hilbert matrix.
However, we will refer to $(H_S(i,j))$ as the Hilbert function.
\end{remark}

We wish to study the Hilbert functions of rings of the form
$R/I$ where $I$ is the ideal associated to a fat point scheme
in $\popo$.  We now recall the relevant definitions.

Let $\pr^1 := \pr^1_{\bf k}$ be the projective line defined over
${\bf k}$, and let $\popo$ be the product space.  The coordinate
ring of $\popo$ is the bigraded ring $R = {\bf
k}[x_0,x_1,y_0,y_1]$ where $\deg x_i = (1,0)$ and $\deg y_i =
(0,1)$.

Suppose that
\[
P = [a_0:a_1] \times [b_0:b_1] \in \popo
\]
is a point in this space.  The ideal $\wp$ associated to $P$ is
the bihomogeneous ideal
\[
\wp = (a_1x_0 - a_0x_1,b_1y_0 - b_0y_1).
\]
The ideal $\wp$ is a prime ideal of height two that is generated
by an element of degree $(1,0)$ and an element of degree $(0,1)$.

If $P = P_1 \times P_2 \in \popo$, then we shall sometimes write
$L_{P_1}$ and $L_{P_2}$ for the generators
of the  ideal $\wp = (L_{P_1},L_{P_2})$ defining $P$ where $L_{P_1}$
is a form of degree $(1,0)$ and $L_{P_2}$ is a form of degree
$(0,1)$.  Since $\popo \cong \mathcal{Q}$, the quadric surface
in $\pr^3$, it is useful to note that $L_{P_1}$ defines a line in
one ruling of $\mathcal{Q}$ and $L_{P_2}$ defines a line in the
other ruling, and $P$ is the point of intersection of these two
lines.

Let $\X$ be a set of $s$ reduced points in $\popo$. Let
$\pi_1:\popo \rightarrow \pr^1$ denote the projection morphism
defined by $P_1 \times P_2 \mapsto P_1$.  Let $\pi_2:\popo
\rightarrow \pr^1$ be the other projection morphism. The set
$\pi_1(\X) = \{R_1,\ldots,R_r\}$ is the set of $r \leq s$  distinct
first coordinates that appear in $\X$.  Similarly, the set
$\pi_2(\X) = \{Q_1,\ldots,Q_t\}$ is the set of $t \leq s$
distinct second coordinates.  For $i = 1,\ldots,r$,
let $L_{R_i}$ denote the $(1,0)$ form that vanishes at all the
points of $\popo$ which have first coordinate $R_i$. Similarly,
for $j = 1,\ldots,t$, let $L_{Q_j}$ denote the $(0,1)$ form that
vanishes at all the points whose second coordinate is $Q_j$.

Let $D:=\{(i,j) ~|~ 1 \leq i \leq r, 1 \leq j \leq t\}.$
If $P \in \X$, then $I_{P} = (L_{R_i},L_{Q_j})$ for some
$(i,j) \in D.$  (Note that this does
not mean that if $(i,j) \in D$, then  $P_{ij} \in \X$.  There may be a
pair $(i,j) \in D$, but $P_{ij} \not\in \X$.)
For each $(i,j) \in D$, let $m_{ij}$ be a
positive integer if $P_{ij} \in \X$, otherwise, let $m_{ij} = 0$.
Then we denote by $Z$ the subscheme of $\popo$ defined by the
saturated bihomogeneous ideal
\[
\Iz = \bigcap_{(i,j) \in D} \wp_{ij}^{m_{ij}}
\]
where $\wp_{ij}^0 := (1)$. We say $Z$ is a {\it fat point
scheme} of $\popo$.  We sometimes say that $Z$ is a {\it set of
fat points}. The integer $m_{ij}$ is called the {\it
multiplicity} of the point $P_{ij}$. We shall sometimes denote
the fat point scheme as
\[
Z = \{ (P_{ij};m_{ij}) ~|~ (i,j) \in D\}.
\]
In the case all the non-zero $m_{ij}$ are the same, we call $Z$ a
{\em{homogeneous fat point scheme}}. The {\it support} of $Z$,
written $\supp(Z)$ is the set of points $\X$. If $\X =
\supp(Z)$, then $\Ix = \sqrt{\Iz}$.

Let $\Iz$ be the defining ideal of a fat point scheme $Z
\subseteq \popo$.  Because the ideal $\Iz \subseteq R$ is a
bihomogeneous ideal we can study its Hilbert function
$H_{R/\Iz}$. We sometimes write $H_{Z}$ to denote $H_{R/\Iz}$,
and say $H_{Z}$ is the {\it Hilbert function of $Z$}.

We give some elementary results about the Hilbert function of a
fat point scheme in $\popo$.  These results generalize some of
the results of ~\cite{VT1} about sets of simple points.

It was shown in \cite[Lemma 3.3]{VT1} that if $\X$ is
a reduced set of points, then there exists a
$(1,0)$ form $L \in R$ (respectively, a $(0,1)$ form $L' \in R$)
that is a non-zero divisor of $R/\Ix$.  The proof of this
lemma can extend to the non-reduced case:

\begin{lemma}   \label{nonzerodivisor}
Let $Z$ be a fat point scheme of $\popo$.  Then there exists a
bihomogeneous element $L \in R$ (respectively, $L' \in R$) with
$\deg L = (1,0)$ (respectively, $\deg L' = (0,1)$) such that
$\ol$ (respectively, $\ol'$) is a non-zero divisor of $R/\Iz$.
\end{lemma}

The existence of these non-zero divisors enables us to prove
the following:
\begin{proposition} \label{elementaryresults}
Let $Z$ be a fat point scheme in $\popo$  and suppose that
$H_{Z}$ is the Hilbert function of $Z$.  Then
\begin{enumerate}
\item[$(i)$]
for all $(i,j) \in \N^2$,
$
H_{Z}(i,j)\leq H_{Z}(i+1,j)$, and
$H_{Z}(i,j)\leq H_{Z}(i,j+1).$
\item[$(ii)$]
if $H_{Z}(i,j) =H_{Z}(i+1,j)$, then $H_{Z}(i+1,j)=H_{Z}(i+2,j)$.
\item[$(iii)$]if $H_{Z}(i,j) = H_{Z}(i,j+1)$, then
$H_{Z}(i,j+1)=H_{Z}(i,j+2)$.
\end{enumerate}
\end{proposition}

\begin{proof}
Let $\ol$ be the non-zero divisor of $R/\Iz$ from Lemma ~\ref{nonzerodivisor}
with $\deg L = (1,0)$.  For any $(i,j) \in \N^2$,
the map
$
(R/\Iz)_{i,j} \stackrel{\times \ol}{\longrightarrow} (R/\Iz)_{i+1,j}
$
is an injective map of vector spaces because $\ol$ is a non-zero
divisor. It then follows that $H_{Z}(i,j) \leq H_{Z}(i+1,j)$ for
all $(i,j) \in \N^2$.  The other statement of $(i)$ is proved
similarly.

The proof of $(ii)$ and $(iii)$ are similar, so we will only show
$(ii)$.  Let $\overline{L}$ be as above.  For each
$(i,j) \in \N^2$, we have
the following short exact sequence of vector spaces:
\[
0 \longrightarrow  \left( R/\Iz\right)_{i,j}
\stackrel{\times \overline{L}}{\longrightarrow} \left( R/\Iz\right)_{i+1,j}
\longrightarrow \left(R/(\Iz,L)\right)_{i+1,j} \longrightarrow 0.
\]
If $H_{Z}(i,j) = H_{Z}(i+1,j)$, then this implies that the
morphism $\times \overline{L}$ is an isomorphism of vector
spaces, and thus, $(R/(\Iz,L))_{i+1,j} = 0$,
or equivalently, $(\Iz,L)_{i+i,j} = R_{i+1,j}$.   But
then $(\Iz,L)_{i+2,j} = R_{1,0} \otimes_{\bf k} R_{i+1,j} = R_{i+2,j}$,
and thus,
$(R/(\Iz,L))_{i+2,j} = 0$ as well.  The
exact sequence then implies that $(R/\Iz)_{i+1,j} \cong (R/\Iz)_{i+2,j}$.
\end{proof}

\begin{remark} Proposition \ref{elementaryresults} implies
that the values in the columns and rows of the Hilbert function
$H_{Z}$, written as a matrix, must eventually stabilize, that is,
stay constant.   However, at least two questions remain.  First,
where do the rows and columns stabilize?  Second, at what values
must the columns and rows stabilize?  These questions are
answered in the following sections (Corollary
\ref{eventualbehaviour}).
\end{remark}

\begin{remark}
Because Lemma \ref{nonzerodivisor} shows the existence of a
non-zero divisor in $R/\Iz$ for any fat point scheme $Z$ of
$\popo$, it follows that the inequality $\operatorname{depth}
R/\Iz \geq 1$ always holds.   It should be noted that the
arguments used in Lemma \ref{nonzerodivisor} and Proposition
\ref{elementaryresults} use nothing special about $\popo$ and
can be extended to fat point schemes in $\pnk$. Proposition
\ref{elementaryresults} could also be deduced from Propositions
2.5 and 2.7 of \cite{GuMaRa1}.
\end{remark}


\section{Fat Point Schemes whose Support is on a Line}

In this section we investigate the Hilbert functions of fat
point schemes in $\popo$ whose support lies on a line defined
either by a form of degree $(1,0)$ or a form of degree $(0,1)$.
Because $\popo \cong \mathcal{Q}$, the quadric surface of
$\pr^3$, this is equivalent to studying those fat point schemes
whose support is on a line of the rulings of the surface. We show
that the Hilbert function in this case can be computed directly
from the multiplicities of the points.
This result is a key component of our proof in the next section
describing the eventual behaviour of all fat point schemes in $\popo$.

So, let $Z$ be the fat point scheme
\[
Z =
\{(P_{11};m_{11}),(P_{12};m_{12}),(P_{13};m_{13}),\ldots,(P_{1s};m_{1s})\}
\]
of $s$ fat points where $P_{1j} = R_1 \times Q_j$. Then
$\supp(Z) = \{P_{11},\ldots,P_{1s}\}$. It follows that
$\supp(Z)$ lies on the line defined by the form $L_{R_1} \in
R_{1,0}$.

Let $Z'$ denote a fat point scheme whose support lies on a line
defined by a form of degree $(0,1)$, that is, $ Z' =\{(Q_1
\times R_1;m_{11}),\ldots,(Q_s \times R_1;m_{s1})\}$ with $Q_i$
and $R_1$ as in $Z$.  Then, for any $(i,j) \in \N^2$,
$(\Iz)_{i,j} \cong (I_{Z'})_{j,i}$, and therefore, $H_{Z}(i,j) =
H_{Z'}(j,i)$. Because of this relation, it is enough to
investigate the case that the support of $Z$ is contained on the
line defined by a form of degree $(1,0)$.

\begin{remark} The following result can be recovered from Theorem
$4.1$ of \cite{GuMaRa1} and Theorem $2.1$ in \cite{Gu}
if one first shows that these schemes are
arithmetically Cohen-Macaulay.  However, we give a new proof of
this result that does not depend on knowing that the scheme is
Cohen-Macaulay.
\end{remark}

\begin{theorem} \label{fatpointsonaline}
Let $Z
=\{(P_{11};m_{11}),(P_{12};m_{12}),\ldots,(P_{1s},m_{1s})\}$ be
a fat point scheme in $\popo$ whose support is on a line defined
by a form of degree $(1,0)$. Set $m = \max\{m_{1j}\}_{j=1}^s$.
For $h = 0,\ldots,m-1$, set $a_h = \sum_{j=1}^s (m_{1j}-h)_+$
where $(n)_+ := \max\{0,n\}$. Then the Hilbert function of $Z$ is
\begin{eqnarray*}
H_{Z} & = &  \bmatrix
1& 2 & \cdots  & a_{0}-1 & a_0 & a_0 & \cdots \\
1& 2 &\cdots  & a_0-1 & a_0 & a_0 & \cdots \\
\vdots & \vdots && \vdots  & \vdots &\vdots&\ddots
\endbmatrix
+ \bmatrix
0& 0 & \cdots  & 0 & 0 & 0 & \cdots\\
1& 2 & \cdots  & a_{1}-1 & a_1 & a_1 & \cdots \\
1& 2 & \cdots  & a_{1}-1 & a_1 & a_1 & \cdots \\
\vdots & \vdots && \vdots  & \vdots &\vdots&\ddots
\endbmatrix \\
& &
+ \cdots +
\bmatrix
0& 0 & \cdots  & 0 & 0 & 0 &\cdots\\
\vdots& \vdots& & \vdots & \vdots & \vdots \\
0& 0 & \cdots  & 0 & 0 & 0 &\cdots\\
1& 2 & \cdots  & a_{m-1}-1 & a_{m-1} & a_{m-1} & \cdots \\
1& 2 & \cdots  & a_{m-1}-1 & a_{m-1} & a_{m-1} & \cdots \\
\vdots & \vdots & &\vdots  & \vdots &\vdots&\ddots
\endbmatrix.
\end{eqnarray*}
\end{theorem}

\begin{proof}
For each $j =1,\ldots,s$, the ideal associated to
$P_{1j}$ is $\wp_{1j} = (L_{R_1},L_{Q_j})$.  Set $L = L_{R_1}$
and note that $L$ defines the $(1,0)$ line
in $\popo$ on which all the points lie.  Now for each
$0 \leq h \leq m-1$ we set
\[
Z_h = \{(P_{11};(m_{11}-h)_+),\ldots,(P_{1s};(m_{1s}-h)_+)\}
\]
and let $I_{Z_h}$ be the associated ideal.  Thus $Z_0 = Z$.  Furthermore,
we have the identity
$L^h \cap \Iz = L^h\cdot I_{Z_h}$
for each $h = 0,\ldots,m-1$.

Since $L^m \in \Iz$, we have
$
0 = \overline{L}^m\cdot S \subseteq
\overline{L}^{m-1}\cdot S \subseteq \cdots \subseteq
\overline{L}\cdot S \subseteq S
$
where $S = R/\Iz$ and $\overline{L}^i$ denotes the image
of $L^i$ in $S$.  It then follows that
\[
H_{Z}(i,j) = \dim_{\bf k} S_{ij} = \sum_{h=0}^{m-1}
\dim_{\bf k} \left(
\frac{\overline{L}^h\cdot S}{\overline{L}^{h+1}\cdot S}\right)_{i,j}.
\]
Now for each $h = 0,\ldots,m-1$,
\[\frac{\overline{L}^h\cdot S}{\overline{L}^{h+1}\cdot S}
\cong
\frac{L^hR}{L^{h+1} + L^h\cap\Iz} \cong \frac{L^hR}{L^{h+1}+ L^hI_{Z_h}}
\cong \overline{L}^h\left(\frac{R}{L+I_{Z_h}}\right).
\]
Hence $\dim_{\bf k} \left(
\frac{\overline{L}^h\cdot S}{\overline{L}^{h+1}\cdot S}\right)_{i,j}
= \dim_{\bf k} \left(R/(L+I_{Z_h})\right)_{i-h,j}$, and
thus
\[H_{Z}(i,j) = \sum_{h=0}^{m-1} \dim_{\bf k}
\left(R/(L+I_{Z_h})\right)_{i-h,j}.\]

To compute $H_Z$, we thus need to compute the Hilbert function
of $R/(L+I_{Z_h})$ for each $h$.  We now note that for
each $h$,
\[
(L+I_{Z_h}) = (L,L_{Q_1}^{(m_{11}-h)_+}\cdots L_{Q_s}^{(m_{1s}-h)_+}),
\]
that is, $(L+I_{Z_h})$ is a complete intersection generated by
forms of degree $(1,0)$ and $(0,a_h)$.  The resolution
of $(L+I_{Z_h})$ is given by the {\it Koszul resolution}, i.e.,
\[
0 \longrightarrow R(-1,-a_h) \longrightarrow R(-1,0)\oplus R(0,-a_h)
\longrightarrow (L+I_{Z_h}) \longrightarrow 0.\]
Hence, the Hilbert function of $R/(L+I_{Z_h})$ is
\[
H_{R/(L+I_{Z_h})}  =   \bmatrix
1& 2 & \cdots  & a_h-1 & a_h & a_h & \cdots \\
1& 2 &\cdots  & a_h-1 & a_h & a_h & \cdots \\
\vdots & \vdots && \vdots  & \vdots &\vdots&\ddots
\endbmatrix.
\]
This now completes the proof.
\end{proof}

From now on, if $\alpha = (a_0,\ldots,a_{m-1})$ is a tuple of
non-negative integers, then by  $a_k \in \alpha$ we shall
mean that $a_k$ appears as a coordinate in $\alpha$.
The following corollary of Theorem \ref{fatpointsonaline} will be required
in the next section.

\begin{corollary}   \label{fatpointsonlinecorollary}
With the notation as in Theorem \ref{fatpointsonaline}, let
$\alpha = (a_0,\ldots,a_{m-1})$.  Fix  $j\in \N$.  Then, for all $
i \geq m-1 = \max\{m_{1k}\}_{k=1}^s -1,$
\begin{eqnarray*}
H_{Z}(i,j) & = &\#\{a_k \in \alpha ~|~ a_k \geq 1\}
+\#\{a_k \in \alpha ~|~ a_k \geq 2\} + \cdots \\
& & +\#\{a_k \in \alpha ~|~ a_k \geq j+1\}.
\end{eqnarray*}
\end{corollary}

\begin{proof}
Fix a $j \in \N$, and set
\[
(*) = \#\{a_k \in \alpha ~|~ a_k \geq 1\}
+\#\{a_k \in \alpha ~|~ a_k \geq 2\} + \cdots +
\#\{a_k \in \alpha ~|~ a_k \geq j+1\}.
\]
From our definition of $a_0,\ldots,a_{m-1}$, it follows
that $a_0 \geq a_1 \geq \cdots \geq a_{m-1}$.  Let $l$ be
the largest index such that $a_0,\ldots,a_{l-1} \geq j+1$
but $a_l,\ldots,a_{m-1} < j+1$.  Set $\alpha' = (a_l,\ldots,a_{m-1}).$

For each integer $h = 1, \ldots, j+1$, we have
\[
\#\{a_k \in \alpha ~|~ a_k \geq h\} = l + \#\{a_k \in \alpha' ~|~ a_k \geq h\}.
\]
Thus
\begin{eqnarray*}
(*) = (j+1)l + \#\{a_k \in \alpha' ~|~ a_k \geq 1\} +
\cdots +  \#\{a_k \in \alpha' ~|~ a_k \geq a_l\}.
\end{eqnarray*}
If we set $(**) =  \#\{a_k \in \alpha' ~|~ a_k \geq 1\} +
\cdots +  \#\{a_i \in \alpha' ~|~ a_k \geq a_l\}$, then
\begin{eqnarray*}
(**)& = &
 \#\{a_k \in \alpha' ~|~ a_k = 1\} +  2 \#\{a_k \in \alpha' ~|~ a_k = 2\} +
\cdots + \\
&& a_l\#\{a_k \in \alpha' ~|~ a_k = a_l\} \\
& =& a_l + a_{l+1} + \cdots + a_{m-1}.
\end{eqnarray*}
Hence, $(*) = (j+1)l + a_l + a_{l+1} + \cdots + a_{m-1}$.

On the other hand, by Theorem ~\ref{fatpointsonaline}, if
$i \geq m-1$, then
$\dim_{\bf k} (R/\Iz)_{i,j} = \sum_{h=1}^{s} \min\{j+1,a_{h}\}$.
Since $a_0,\ldots,a_{l-1} \geq j+1$, it follows that
\[
\dim_{\bf k} (R/\Iz)_{i,j} = (j+1)l + a_l + a_{l+1} + \cdots a_{m-1} = (*)
\]
which is what we wished to prove.
\end{proof}


\section{The Eventual Behaviour of the Hilbert Function of a Fat Point Scheme}

Let $P_1,\ldots,P_s$ be $s$ distinct points of $\popo$ and
suppose $m_1,\ldots,m_s$ are arbitrary positive integers. Let $Z
= \{(P_1;m_1),\ldots,(P_s;m_s)\}$ be the resulting fat point
scheme of $\popo$. In this section we wish to describe the
eventual behaviour of the Hilbert function of $Z$. We will show
that the eventual values of the Hilbert function depend only
upon the numbers $m_1,\ldots,m_s$ and numerical information
describing $\X = \supp(Z)$. This result is a generalization of a
result of the second author \cite[Corollary 5.13]{VT1} about
sets of points in $\popo$.

We start by defining our notation.
If $Z$ is a fat point scheme, let $\X$ denote the support
of $Z$.  We suppose that $|\X|=s$. Let $\pi_1(\X)$ and $\pi_2(\X)$ be
defined as in the previous section.  For each $R_i \in \pi_1(\X)$,
define
\[
Z_{1,R_i} :=
\{(P_{ij_{1}};m_{ij_{1}}),(P_{ij_{2}};m_{ij_{2}}),\ldots,
(P_{ij_{\alpha_i}};m_{ij_{\alpha_i}})\}
\]
where $P_{ij_k} = R_i \times Q_{j_k}$
are those points of $\supp(Z)$ whose first projection is
$R_i$.  Thus $\pi_1(\supp(Z_{1,R_i})) = \{R_i\}$, and furthermore
it follows that
\[\Iz = \bigcap_{i=1}^r I_{Z_{1,R_i}}.\]
For each $R_i \in \pi_1(\X)$
define $l_i := \max\{m_{ij_1},\ldots,m_{ij_{\alpha_i}}\}$.  Then,
for each integer $0 \leq k \leq l_i-1$, we define
\[
a_{i,k} := \sum_{j=1}^{\alpha_i} (m_{ij} - k)_+ \hspace{.5cm}
\mbox{where $(n)_+ := \max\{n,0\}$.}
\]
Let $\alpha_{R_i} := (a_{i,0},\ldots,a_{i,l_i-1})$ for
each $R_i \in \pi_1(\X)$.  Define
\begin{eqnarray*}
\az & := &(\alpha_{R_1},\ldots,\alpha_{R_r}) \\
& =&
(a_{1,0},\ldots,a_{1,l_1-1},a_{2,0},\ldots,a_{2,l_2-1},\ldots,
a_{r,0},\ldots,a_{r,l_r-1}).
\end{eqnarray*}

Similarly, for each $Q_j \in \pi_2(\X)$,
define
\[
Z_{2,Q_j} :=
\{(P_{i_{1}j};m_{i_{1}j}),(P_{i_{2}j};m_{i_{2}j}),\ldots,
(P_{i_{\beta_j}j};m_{i_{\beta_j}j})\}
\]
where $P_{i_kj} = R_{i_k} \times Q_j$ are those points of
$\supp(Z)$ whose second projection is $Q_j$.
Thus $\pi_2(\supp(Z_{2,Q_j})) = \{Q_j\}$.
For $Q_j \in \pi_2(\X)$
define $l'_j = \max\{m_{i_1j},\ldots,m_{i_{\beta_j}j}\}$.  Then,
for each integer $0 \leq k \leq l'_j-1$, we define
\[
b_{j,k} := \sum_{i=1}^{\beta_j} (m_{ij} - k)_+ \hspace{.5cm}
\mbox{where $(n)_+ := \max\{n,0\}$.}
\]
Let $\beta_{Q_j} := (b_{j,0},\ldots,b_{j,l'_j-1})$ for each $Q_j \in
\pi_2(\X)$.  Define
\begin{eqnarray*}
\bz & := &(\beta_{Q_1},\ldots,\beta_{Q_t}) \\
& =& (b_{1,0},\ldots,b_{1,l'_1-1},b_{2,0},\ldots,b_{2,l'_2-1},\ldots,
b_{t,0},\ldots,b_{t,l'_t-1}).
\end{eqnarray*}

\begin{example} \label{multipointexample}
With the above notation, let us determine the tuples $\az$ and $\bz$ associated
to the
scheme $Z =
\{(P_{11};4),(P_{12};2),(P_{23};3),(P_{32};2),(P_{41};3)\}$.
The subscheme $Z_{1,R_1}$ is
\[
Z_{1,R_1} =\{ (P_{11};4),(P_{12};2)\}.
\]
We set $l_1 := \max \{4,2\} = 4$.  Then
\begin{eqnarray*}
a_{1,0} & = & 4 + 2 = 6 \\
a_{1,1} & = & (4-1)_+ + (2-1)_+ = 4\\
a_{1,2} & = & (4-2)_+ + (2-2)_+ = 2\\
a_{1,3} & = & (4-3)_+ + (2-3)_+ = 1.
\end{eqnarray*}
Hence, $\alpha_{R_1} = (6,4,2,1)$.
For $R_2,R_3,$ and $R_4$, we get
$\alpha_{R_2} =
(3,2,1)$, $\alpha_{R_3} = (2,1)$, $\alpha_{R_4} = (3,2,1)$. Hence
\[
\az = (6,4,2,1,3,2,1,2,1,3,2,1).
\]
Similarly, for $Q_1,Q_2,Q_3 \in \pi_2(\X)$, $l'_1 = 4$, $l'_2 =
2$ and $l'_3 = 3$. So, we have $\beta_{Q_1} = (7,5,3,1)$,
$\beta_{Q_2} = (4,2)$, and $\beta_{Q_3} = (3,2,1)$, and
therefore,
\[
\bz = (7,5,3,1,4,2,3,2,1).
\]
\end{example}

We now state and prove our main result about the eventual behaviour of the
Hilbert function.  Recall that if we write $a_k \in \alpha$,
where $\alpha$ is a tuple of non-negative integers, then we shall
mean that $a_k$ appears as a coordinate in $\alpha$.

\begin{theorem} \label{borderresult}
Let $Z$ be a fat point scheme of $\popo$.
Then, with the above notation,
\begin{enumerate}
\item[$(i)$]
for a fixed $j\in \N$, if $i \geq (l_1 + \cdots +l_r) -1,$ then
\begin{eqnarray*}
\dim_{\bf k} (R/\Iz)_{i,j} & =&
\#\{a_{k,l} \in \az ~|~ a_{k,l} \geq 1\} +
\#\{a_{k,l} \in \az ~|~ a_{k,l} \geq 2\} + \cdots \\
&&+\#\{a_{k,l} \in \az ~|~ a_{k,l} \geq j+1\}.
\end{eqnarray*}

\item[$(ii)$]for a fixed $i\in \N$, if $j \geq (l'_1 + \cdots + l'_t) -1,$ then
\begin{eqnarray*}
\dim_{\bf k} (R/\Iz)_{i,j} & =&
\#\{b_{k,l} \in \bz ~|~ b_{k,l} \geq 1\} +
\#\{b_{k,l} \in \bz ~|~ b_{k,l} \geq 2\} + \cdots \\
&&+\#\{b_{k,l} \in \bz ~|~ b_{k,l} \geq i+1\}.
\end{eqnarray*}
\end{enumerate}
\end{theorem}

\begin{proof}
We will only prove $(i)$ since the proof of statement of $(ii)$
is similar.  Let $Z$ be a set of fat points in $\popo$, and let
$\X = \supp(Z)$.  The proof is by induction on $r = |\pi_1(\X)|$.
If $r = 1$, i.e.,
$\pi_1(\X) = \{R_1\}$, the conclusion follows from Corollary
\ref{fatpointsonlinecorollary}.

So, suppose that $r > 1$, and the theorem holds for all fat point
schemes $Z'$ with $|\pi_1(\supp(Z'))| < r$.  For each $R_i \in
\pi_1(\X)$, we let $I_{Z_{1,R_i}}$ denote the ideal that defines
the subscheme $ Z_{1,R_i} :=
\{(P_{ij_{1}};m_{ij_{1}}),(P_{ij_{2}};m_{ij_{2}}),\cdots,
(P_{ij_{\alpha_i}};m_{ij_{\alpha_i}})\}. $ We set
\[
I_{\Y_1} := \bigcap_{i=1}^{r-1} I_{Z_{1,R_i}}
 \hspace{.5cm}\mbox{and}
\hspace{.5cm} I_{\Y_2} :=  I_{Z_{1,R_r}}.
\]
The ideals $I_{\Y_1}$ and $I_{\Y_2}$ are the defining ideals of
fat point schemes in $\popo$ with $|\pi_1(\supp(\Y_i))| < r$ for
$i = 1,2$.  We shall also require the following result about
$I_{\Y_1}$ and $I_{\Y_2}$.

\noindent {\it Claim. }  For any $j \in \N$, if $i \geq l_1 +
\cdots + l_r -1$, then $(I_{\Y_1} + I_{\Y_2})_{i,j} = R_{i,j}$.

\noindent {\it Proof of the Claim. } Set $m = l_1 + \cdots +
l_r$.  It is enough to show that $(I_{\Y_1} + I_{\Y_2})_{m-1,0} =
R_{m-1,0}$.
Recall that for each $R_i \in \pi_1(\X)$, the integer $l_i$ is
defined to be $l_i = \max\{m_{ij_c}\}_{c=1}^{\alpha_i}$ where
$Z_{1,R_i}$ is as above. If $(L_{R_i},L_{Q_{j_c}})$ is the ideal
associated to the point $P_{ij_c}$, then $I_{Z_{1,R_i}} =
\bigcap_{c=1}^{\alpha_i} (L_{R_i},L_{Q_{j_c}})^{m_{ij_c}}$. Note
that $\deg L_{R_i} = (1,0)$ and $\deg L_{Q_{j_c}} = (0,1)$. From
this description of $I_{Z_{1,R_i}} $, it follows that
$L_{R_i}^{l_i} \in I_{Z_{1,R_i}}$.  Thus $L_{R_1}^{l_1}\cdots
L_{R_{r-1}}^{l_{r-1}} \in I_{\Y_1}$ and $L_{R_r}^{l_r} \in
I_{\Y_2}$.

Set $J: = (L_{R_1}^{l_1}\cdots
L_{R_{r-1}}^{l_{r-1}},L_{R_r}^{l_r}) \subseteq I_{\Y_1} +
I_{\Y_2}$.  Since $J$ is generated by a regular sequence, the
bigraded resolution of $J$ is given by the Koszul resolution:
\[
0 \longrightarrow R(-m,0) \longrightarrow R(-m+l_r,0) \oplus
R(-l_r,0) \longrightarrow J \longrightarrow 0.
\]
If we use this exact sequence to calculate the dimension of $J_{m-1,0}$,
then we find
\begin{eqnarray*}
\dim_{\bf k} J_{m-1,0} & = & (m-1 -(m-l_r) +1) + (m-1 - l_r +1)
- (m-1 -m +1) \\
& = & l_r + m - l_r = m = \dim_{\bf k} R_{m-1,0}.
\end{eqnarray*}
Since $\dim_{\bf k} J_{m-1,0} \leq
\dim_{\bf k} (I_{\Y_1} + I_{\Y_2})_{m-1,0} \leq \dim_{\bf k} R_{m-1,0}$,
the conclusion $ (I_{\Y_1} + I_{\Y_2})_{m-1,0} = \dim_{\bf k} R_{m-1,0}$
now follows.\hfill$\diamond$

\noindent
From the short exact sequence
\[
0 \longrightarrow I_{\Y_1} \cap I_{\Y_2} = \Iz \longrightarrow
 I_{\Y_1} \oplus I_{\Y_2}\longrightarrow I_{\Y_1}+ I_{\Y_2} \longrightarrow
0
\]
we deduce that
\[
\dim_{\bf k} (\Iz)_{i,j} = \dim_{\bf k} (I_{\Y_1})_{i,j}
+ \dim_{\bf k} (I_{\Y_2})_{i,j} -   \dim_{\bf k} (I_{\Y_1} +
I_{\Y_2})_{i,j}
\]
for all $(i,j) \in \N^2$.  Thus, if $i \geq l_1+\cdots+ l_r -1$, then by
the claim we have
\begin{eqnarray*}
H_{Z}(i,j) & =& (i+1)(j+1) - \dim_{\bf k} (I_{\Y_1})_{i,j} -
\dim_{\bf k} (I_{\Y_2})_{i,j} + \dim_{\bf k} (I_{\Y_1} + I_{\Y_2})_{i,j} \\
& = & (i+1)(j+1) - \dim_{\bf k} (I_{\Y_1})_{i,j} +
(i+1)(j+1) -\dim_{\bf k} (I_{\Y_2})_{i,j} \\
& = & H_{\Y_1}(i,j) + H_{\Y_2}(i,j).
\end{eqnarray*}

For each $h = 1,\ldots,j+1$, it follows that
\[
\#\{a_{k,l} \in \az ~|~ a_{k,l} \geq h \} = \#\{a_{k,l} \in \alpha_{\Y_1} ~|~
a_{k,l} \geq h \} + \#\{a_{t,l} \in \alpha_{\Y_2} ~|~
a_{t,l} \geq h \}
\]
where $\alpha_{\Y_i}$ is the tuple associated to the fat point
scheme $\Y_i$ for $i = 1,2$. The conclusion now follows by the
induction hypothesis and the fact that $H_{Z}(i,j) =
H_{\Y_1}(i,j) + H_{\Y_2}(i,j)$ if $i \geq l_1+\cdots+ l_r -1$.
\end{proof}

\begin{remark}
Suppose that $Z$ is a set of simple points in $\popo$, i.e.,
the multiplicity of each point in $Z$ is one.  So, if $\pi_1(Z)
= \{R_1,\ldots,R_r\}$, then $Z_{1,R_i} = \{ R_i \times
Q_{i_1},\ldots,R_i \times Q_{i_{\alpha_i}}\}$ for $i =
1,\ldots,r$. So, $l_i =1$, and thus, $a_{i,0} =
\sum_{j=1}^{\alpha_i} 1 = \alpha_i$. So, $\az =
(\alpha_1,\ldots,\alpha_r)$, which is exactly how $\az$ is
defined for sets of simple points in ~\cite{VT1}.
Thus Theorem ~\ref{borderresult} generalizes
\cite[Proposition 5.11]{VT1} for sets of points in $\popo$
to fat point schemes in $\popo$.
\end{remark}

We can rewrite
Theorem \ref{borderresult} more succinctly.

\begin{corollary}   \label{eventualbehaviour}
Let $Z$ be a fat point scheme in $\popo$. With the notation as
in  Theorem ~\ref{borderresult}, let $m = l_1 + \cdots + l_r $
and  $m' = l'_1 + \cdots + l'_t$. Then
\[
H_{Z}(i,j) = \left\{
\begin{array}{ll}
\sum_{i=1}^s \binom{m_i +1}{2} & \mbox{if $(i,j) \geq (m-1,m'-1)$} \\
H_{Z}(m-1,j) & \mbox{if $i \geq m-1$ and $j < m'-1$} \\
H_{Z}(i,m'-1) & \mbox{if $j \geq m'-1$ and $i < m-1$}
\end{array}
\right..
\]
\end{corollary}

\begin{proof}
For any $j \in \N$, if $i \geq m-1$, then Theorem \ref{borderresult}
implies that $H_{Z}(i,j) = H_{Z}(m-1,j)$.  Similarly,
for any $i \in \N$, if $j \geq m'-1$, then $H_{Z}(i,j) = H_{Z}(i,m'-1)$.
Thus, for any $(i,j) \geq (m-1,m'-1)$,
we  have $H_{Z}(i,j) = H_{Z}(i,m'-1)
= H_{Z}(m-1,m'-1)$.

All that remains to be shown is that
$H_{Z}(m-1,m'-1) = \sum_{i=1}^s \binom{m_i +1}{2}$.
From Theorem ~\ref{borderresult} it follows that
\begin{eqnarray*}
H_{Z}(m-1,j) & = &\#\{a_{k,l} \in \az ~|~ a_{k,l} \geq 1\} +
\cdots
+\#\{a_{k,l} \in \az ~|~ a_{k,l} \geq j+1\} \\
& = & \#\{a_{k,l} \in \az ~|~ a_{k,l} = 1\} +
2\#\{a_{k,l} \in \az ~|~ a_{k,l} = 2\} +\cdots + \\
& & (j+1)\#\{a_{k,l} \in \az ~|~ a_{k,l} = j+1\}.
\end{eqnarray*}
Thus, if $j \gg 0$, then $ H_{Z}(m-1,j) =
\sum_{k=1}^r\sum_{l=1}^{l_k-1} a_{k,l}. $ For any  $k \in
\{1,\ldots,r\}$
\begin{eqnarray*}
\sum_{l=1}^{l_k-1} a_{k,l} & = & a_{k,0} + a_{k,1} + \cdots + a_{k,l_k-1} \\
& = & \left[m_{i_1} + (m_{i_1}-1) + \cdots + 2 + 1 \right]
+ \cdots +  \left[m_{i_{\alpha_i}} + (m_{i_{\alpha_i}}-1) +
\cdots + 2 + 1 \right] \\
& = & \binom{m_{i_1} + 1}{2} + \cdots +\binom{m_{i_{\alpha_i}} + 1}{2}.
\end{eqnarray*}
It then follows that $H_{Z}(m-1,j) = \sum_{i=1}^s \binom{m_i
+1}{2}$ if $j \gg 0$.  In particular, $H_{Z}(m-1,m'-1) =
\sum_{i=1}^s \binom{m_i
+1}{2}$.
\end{proof}

\begin{remark}
From the above corollary, we see that if we know the values of
$H_{Z}(m-1,j)$ for $j=0,\ldots,m'$ and the values of
$H_{Z}(i,m'-1)$ for $i = 0,\ldots,m$, then we know the entire
Hilbert function except at a finite number of values.  This
observation motivates the next definition.
\end{remark}

\begin{definition}  Let $Z$ be a fat point scheme and let $\az$ and $\bz$
be constructed as described above.  If $m = |\az|$ and $m' = |\bz|$, then
define the following tuples:
\[
B_C = (H_{Z}(m-1,0),H_{Z}(m-1,1),\ldots,H_{Z}(m-1,m'-1))
\]
and
\[
B_R = (H_{Z}(0,m'-1),H_{Z}(1,m'-1),\ldots,H_{Z}(m-1,m'-1)).
\]
The tuple $B_C$ is called the {\it eventual column vector}
because it contains the values at which the columns will
stabilize.  Similarly, $B_R$ is the {\it eventual row vector}.
Set $B_{Z} := (B_C,B_R)$.  The tuple $B_{Z}$ is called the {\it
border} of the Hilbert function of $Z$.
\end{definition}

The notion of a border was first introduced in
\cite{VT1} for sets of simple points in $\pnk$.  The name is
used to describe the fact that once we know the values of border, then
we know all the values of the Hilbert function ``outside''
the border.  Thus only values ``inside'' the border, i.e., those
$(i,j) \in \N^2$ with $(i,j) \leq (m-1,m'-1)$, need to be calculated
to completely determine the entire Hilbert function.

It follows from Theorem \ref{borderresult} that
the border can be computed directly from the tuples $\az$ and $\bz$.
By borrowing some terminology from combinatorics, we can
make this connection explicit.  Our main reference
for this material is Ryser \cite{R}.  But first, for the remainder
of this paper, we will adopt the following convention about $\az$
and $\bz$.

\begin{convention}
Let $Z$ be a fat point scheme in $\popo$, and suppose that $\az$
and $\bz$ are constructed from $Z$ as described above.  We will
assume that the entries of $\az = (\alpha_1,\ldots,\alpha_m)$
have been reordered so that $\alpha_i \geq \alpha_{i+1}$ for
each $i$.  We assume the same for $\bz$.
\end{convention}

\begin{definition}
\label{partitiondefinition}
A tuple $\lambda = (\lambda_1,\ldots,\lambda_r)$ of positive integers
is a {\it partition}
of an integer $s$ if $\sum \lambda_i = s$ and $\lambda_i \geq
\lambda_{i+1}$ for every $i$.  We write $\lambda =
(\lambda_1,\ldots,\lambda_r) \vdash s$.  The {\it conjugate} of $\lambda$
is the tuple $\lambda^* = (\lambda^*_1,\ldots,\lambda^*_{\lambda_1})$
where $\lambda_i^* = \#\{\lambda_j \in \lambda ~|~ \lambda_j \geq i\}$.
Furthermore, $\lambda^* \vdash s$.
\end{definition}

\begin{example}
If $Z = \{(P_1,m_1),\ldots,(P_s,m_s)\}$ is a fat point scheme of
$\popo$, then the tuples $\az$ and $\bz$ are partitions of $\deg
Z = \sum_{i=1}^s \binom{m_i +1}{s}$.
\end{example}

\begin{definition}
\label{ferrers}
To any partition $\lambda = (\lambda_1,\ldots,\lambda_r) \vdash s$
we can associate the following diagram:  on an $r \times \lambda_1$
grid, place $\lambda_1$ points on the first line,
$\lambda_2$ points on the second, and so on.  The resulting diagram is
called the {\it Ferrer's diagram} of $\lambda$.
\end{definition}

\begin{example}
Suppose $\lambda = (4,4,3,1) \vdash 12$.  Then the Ferrer's diagram
is
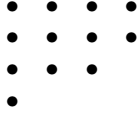
\[
\begin{tabular}{cccc}
$\bullet$ & $\bullet$ & $\bullet$ & $\bullet$ \\
$\bullet$ & $\bullet$ & $\bullet$ & $\bullet$ \\
$\bullet$ & $\bullet$ & $\bullet$ & \\
$\bullet$ & & &
\end{tabular}\]
The conjugate of $\lambda$ can be read off the Ferrer's diagram by
counting the number of dots in each column as opposed to each row.  In
this example $\lambda^* = (4,3,3,2)$.
\end{example}

For any tuple $p: = (p_1,\ldots,p_k)$, we define $\Delta p : = (p_1,
p_2-p_1,\ldots,p_k - p_{k-1})$.

\begin{corollary}   \label{bordercor}
Let $Z$ be a fat point scheme of $\popo$. Then
\begin{enumerate}
\item[$(i)$] $\Delta B_C = \az^*$.
\item[$(ii)$] $\Delta B_R = \bz^*$.
\end{enumerate}
\end{corollary}

\begin{proof}
We use Theorem \ref{borderresult} to calculate $\Delta B_C$:
\[
\Delta B_C = (\#\{\alpha_i \in \az ~|~ \alpha_i \geq 1\},
\#\{\alpha_i \in \az ~|~ \alpha_i \geq 2\},\ldots, \#\{\alpha_i
\in \az ~|~ \alpha_i \geq m'\})
\]
where $m' = |\bz|$.  Since $\#\{\alpha_i \in \az ~|~ \alpha_i
\geq h\}$ is by definition the $h^{th}$ coordinate of
$\az^*$, we have $\Delta B_C = \az^*.$  The proof of $(ii)$ is
the same.
\end{proof}

\begin{remark}
Corollary ~\ref{bordercor} implies that we can compute the Hilbert
function of $Z$ at all but a finite number of values from
only the multiplicities and the relative positions of the points.
\end{remark}

\begin{example}
This example illustrates that in $\popo$ subschemes with the
same border can have different Hilbert functions. Set $R_i = Q_i
= [1:i] \in \pr^1$, and let $P_{ij}$ denote the point $R_i
\times Q_j$. Let
\begin{eqnarray*}
Y_1 & = &
\{(P_{11};1),(P_{22};1),(P_{33};1),(P_{45};1)\} ~\mbox{and}\\
Y_2 & = &\{(P_{11};1),(P_{22};1),(P_{33};1),(P_{44},1) \}.
\end{eqnarray*}
As an exercise one can verify that $\alpha_{Y_1} = \alpha_{Y_2} =
(1,1,1,1)$ and $\beta_{Y_1} = \beta_{Y_2} =
(1,1,1,1)$.  Thus, the two schemes have the same
border.  The Hilbert function of $H_{Y_1}$
is
\[
\bmatrix
1  & 2   & 3  & 4  & 4  & \cdots \\
2  & 4   & 4  & 4  & 4& \cdots \\
3  & 4   & 4  & 4 & 4 & \cdots \\
4  & 4   & 4 &  4 & 4 & \cdots \\
4  & 4   & 4 &  4 & 4 & \cdots \\
\vdots&\vdots & \vdots & \vdots & \vdots & \ddots \\
\endbmatrix
\]
from which we deduce that $(I_{Y_1})_{1,1} = 0$.  On the
other hand, the unique $(1,1)$-form $(x_0y_1 - y_0x_1)$ which passes
through $P_{11},P_{22},$ and $P_{33}$ also passes
through the point $P_{44}$ but not $P_{45}$.  Thus $(I_{Y_2})_{1,1} \neq 0$,
and hence, $H_{Y_1} \neq H_{Y_2}$.
\end{example}

As we have seen, the tuples $\az$ and $\bz$ give us a lot of
information about the Hilbert function of $Z$.  It is therefore
natural to ask which tuples can arise from a fat point scheme $Z$
in $\popo$. Because of Corollary \ref{bordercor}, this is
equivalent to asking what can be the border of the Hilbert
function of a fat point scheme in $\popo$.  The following theorem
places a necessary condition on the tuples $\az$ and $\bz$. We
require the following definition.

\begin{definition}
\label{majorizes}
Let $\lambda = (\lambda_1,\ldots,\lambda_t)$ and $\delta =
(\delta_1,\ldots,\delta_r)$ be two partitions of $s$.  If one
partition is longer, we add zeroes to the shorter one until they
have the same length.  We say $\lambda$ {\it majorizes} $\delta$,
written $\lambda \unrhd \delta$, if
\[ \lambda_1 + \cdots + \lambda_i \geq \delta_1 + \cdots + \delta_i
\mbox{ for $i = 1, \ldots,\max\{t,r\}$}.\]
Majorization induces a partial ordering on the set
of all partitions of $s$.
\end{definition}

\begin{theorem}
Let $Z$ be a scheme of fat points in $\popo$.  Then
\[
\az^* \unrhd \bz.
\]
\end{theorem}
\begin{proof}
We work by induction on $m = |\alpha_{Z}|$.  If $m = 1$, then
$Z$ is a scheme of simple points in $\popo$.  Thus $\az^* \unrhd
\bz$ by Theorem 5.16 in \cite{VT1}.

So, let us suppose that $m > 1$.  We can write $Z$ as
\[
Z = \{(P_{ij};m_{ij}) ~|~ 1 \leq i \leq r, ~1\leq j \leq t \}
\]
where $m_{ij} \geq 0$ and $P_{ij} = R_i \times Q_j$ for some
$R_i, Q_j \in \pr^1$. Recall that if $m_{ij} = 0$, then $P_{ij}
\not\in \supp(Z)$.

For each $i = 1,\ldots,r$, set $m_i := \sum_{j=1}^t m_{ij}$.
After relabeling the $P_{ij}$'s, we can assume that $m_1 =
\max\{m_1,\ldots,m_r\}$. Furthermore, we can also suppose that
after relabeling, $m_{1j}\neq 0$ for $j = 1, \ldots,k$, and
$m_{1j} = 0$ for $j = k+1,\ldots,t$.  Thus $m_1 = m_{11} +
\cdots +m_{ik}$. Note that $m_1 = \alpha_1$, the first
coordinate of $\az$.

Let $\Y$ be the following subscheme of $Z$:
\[
\Y := \{(P_{ij};m'_{ij}) ~|~ 1 \leq i \leq r, ~1\leq j \leq t \}
\]
where
\[
m'_{ij} = \left\{
\begin{array}{ll}
(m_{ij}-1)_+ & i = 1, ~1 \leq j \leq t\\
m_{ij} & 2 \leq i \leq r, ~1 \leq j \leq t
\end{array}
\right.
\]
with $(n)_+ := \max\{0,n\}.$  The subscheme $\Y$ is constructed
from $Z$ by subtracting 1 from the multiplicity of each point on
the $(1,0)$ line that corresponds to $\alpha_1$ in $\az$.

Since $\az = (\alpha_1,\ldots,\alpha_m)$, and because $\alpha_1
= m_1$, from our construction of $\Y$ it follows that
$\alpha_{\Y} = (\alpha_2,\ldots,\alpha_m)$. Therefore, by
induction $\alpha_{\Y}^* \unrhd \beta_{\Y}$.

Let $\beta_{\Y}$ and $\beta_{Z}$ be the tuples associated to
$\Y$ and $Z$, respectively, but for the moment we assume that
$\beta_{\Y}$ and $\beta_{Z}$ have been constructed as first
described at the beginning of Section 3, that is, $\beta_{\Y}$
and $\beta_{Z}$ have not been ordered.

We now describe how  $\bz$ and $\beta_{\Y}$ are related.
Suppose $\beta_{Z} = (b_1,b_2,\ldots,b_l)$ and $\beta_{\Y} =
(b'_1,b'_2,\ldots,b'_h)$.  Clearly $h \leq l$.

If $h=l$, then
\[
b_p=b'_p+1\quad\text{for all $p=1,\ldots,l$.}
\]

If $h < l$, we first insert $(l-h)$ zeroes into the tuple
$\beta_{\Y}$ at specific locations. For $j = 1,\ldots,t$, set
$l'_j := \max\{m_{1j},m_{2j},\ldots,m_{rj}\}$, and for
$d=1,\ldots,t$, set $h_d := \sum_{s=1}^d l'_s$.  Then we insert
a zero into the $h_d^{th}$ spot of $\beta_{\Y}$ if $l'_d =
m_{1d}$ but $l'_d > m_{id}$ for all $i=2,\ldots,r$. It then
follows from our definition of $\Y$ that we are only adding
$(l-h)$ zeroes to $\beta_{\Y}$. Relabel our tuple as $\beta_{\Y}
= (c_1,\ldots,c_l)$.

From our construction of $\Y$ from the scheme $Z$, it follows
that
\[b_i =  \begin{cases}
c_i + 1 &\text{ for } i = 1,\ldots,m_{11}, l'_1+1,\ldots,m_{12},\\
{}& l'_1 + l'_2 + 1,\ldots,m_{13},~\ldots, l'_1
+ l'_2 + \cdots + l'_{k-1}+ 1,\ldots,m_{1k}\\
c_i & \mbox{otherwise}
\end{cases}
\]
So $\beta_{Z}$ can be constructed from $\beta_{\Y}$ by adding 1
to  $m_{11} + m_{12} + \cdots + m_{1k} = m_1 = \alpha_1$
distinct coordinates in $\bz$, and then reordering so that
$\beta_{Z}$ is a partition.

Since $\az = (\alpha_1,\ldots,\alpha_m)$ and $\alpha_{\Y} =
(\alpha_2,\ldots,\alpha_m)$, $\az^*$ can be computed from
$\alpha_{\Y}^*$ by adding 1 to the first $\alpha_1$ entries of
$\alpha_{\Y}^*$.  (If $|\alpha_{\Y}^*| < \alpha_1$, we extend
$\alpha_{\Y}^*$ by adding zeroes so $|\alpha_{\Y}^*| =
\alpha_1$.) By induction, $\alpha_{\Y}^* \unrhd \beta_{\Y}$.  So,
if $\beta_{\Y} = (c_1,\ldots,c_l)$, then
\[
\az^* \unrhd (c_1 + 1,\ldots,c_{\alpha_1} + 1, c_{\alpha_1
+1},\ldots,c_l).
\]
But since $\bz$ can be recovered from $\beta_{\Y}$ by adding 1
to $m_1 = \alpha_1$ distinct entries of $\beta_{\Y}$ (and not
necessarily the first $\alpha_1$ entries) and then reordering,
we have
\[
\az^* \unrhd (c_1 + 1,\ldots,c_{\alpha_1} + 1, c_{\alpha_1
+1},\ldots,c_l) \unrhd \bz.
\]
Hence $\az^* \unrhd \bz$, as desired.
\end{proof}


\section{ACM Fat Point Schemes}

For any fat point scheme in $\pr^n$, the associated coordinate
ring is always Cohen-Macaulay.  In contrast, fat point schemes
in $\pnk$ with $k \geq 2$ may fail to have this property, even
if the support is ACM. See \cite{GuMaRa1,Gu,VT2}
for more details on ACM zero-dimensional schemes in $\pnk$.

A fat point scheme is said to be {\it arithmetically
Cohen-Macaulay} (ACM for short) if the associated coordinate
ring is Cohen-Macaulay. ACM schemes on a smooth quadric
$\mathcal{Q} \cong \popo$ were studied in \cite{GuMaRa1} and by
the first author in \cite{Gu} (which is based on \cite{Gu2}). In
\cite{GuMaRa1} the authors gave a characterization of ACM
schemes in terms of their Hilbert functions. In \cite{Gu}, ACM
fat points schemes in $\popo$ were characterized in terms of the
multiplicities of the points. In this section we show that ACM
schemes can also be classified using the tuples $\az$ and $\bz$
introduced in the previous section.  We will also show how these
various classifications are related.

We begin by recalling the construction and main result of
\cite{Gu}.  Let $Z$ be a fat point scheme in $\popo$ where $Z =
\{(P_{ij};m_{ij}) ~|~ 1 \leq i \leq r, ~1 \leq j \leq t\}$ with
$m_{ij} \geq 0$ and $P_{ij} = R_i \times Q_j$ for some $R_i, Q_j
\in \pr^1$. For each $h \in \N$, and for each tuple $(i,j)$ with
$1 \leq i \leq r$ and $1 \leq j \leq t$, define
\[
t_{ij}(h) := (m_{ij}-h)_+ = \max\{0,m_{ij}-h\}.
\]
The set $\Ssz$ is then defined to be the set of $t$-tuples
\[
\Ssz=\{(t_{i1}(h),\ldots,t_{it}(h)) ~|~  1\leq i \leq r, ~h \in \N\}
\]
For each integer $1 \leq i \leq r$, set $l_{i} :=
\max\{m_{i1},\ldots,m_{it}\}$. For any fat point scheme, we then
have $|\Ssz|= m : =\sum_{i=1}^r l_i$. For each $i = 1,\ldots,r$
and for all $h \in \N$ we set
\[
z_{i,h} :=\sum_{j=1}^{t}{t_{ij}(h)}.
\]
We then define $u_1 := \max_{i,h}\{z_{i,h}\}$, and we recursively define
\[
u_p:=\max_{{i,h}}\{\{z_{i,h}\}\setminus \{u_1,\dots,u_{p-1}\}\} \quad
\text{for}\quad p=2,\dots,m.
\]

\begin{definition}
Let $H_{Z}:\N^2 \rightarrow \N$ be the Hilbert function of a fat
point scheme $Z$ in $\popo$.  The {\it first difference
function} of $H_{Z}$, denoted $\Delta H_{Z}$, is the function
defined by
\[
\Delta H_{Z}(i,j) = H_{Z}(i,j) - H_{Z}(i-1,j) - H_{Z}(i,j-1) +
H_{Z}(i-1,j-1)
\]
where $H_{Z}(i,j) = 0$ if $(i,j) \not\geq (0,0)$.
\end{definition}

With this notation we can state the main result of \cite{Gu}.

\begin{theorem}[{\cite[Theorem 2.1]{Gu}}]   \label{acm}
Let $Z$ be a fat point scheme on $\mathcal{Q} \cong \popo$. Then
the set $\mathcal S_{Z}$ is totally ordered if and only if $Z$ is
ACM. In this case, the first difference function of $H_{Z}$ is:
\[
\Delta H_{Z} = \bmatrix
\underbrace{1 ~ 1 ~ 1 ~ \cdots ~ 1}_{u_1} & 0 & \cdots \\
\underbrace{1 ~ 1 ~ \cdots ~ 1}_{u_2} ~ 0 & 0 & \cdots \\
 \vdots \\
\underbrace{1 ~ \cdots ~ 1}_{u_m} ~ 0 ~ 0 & 0 & \cdots \\
0 ~ \cdots ~ 0 ~0 ~ 0 & 0 & \cdots \\
\vdots  \hspace{.8cm}\vdots ~~ \vdots ~~ \vdots & \vdots & \ddots \\
\endbmatrix
\]
where $u_1,\ldots,u_m$ are defined as above.
\end{theorem}

\begin{remark} \label{az=sz}
From the construction of $u_1,\ldots,u_m$, one can verify that the
identity $\az = (u_1,\ldots,u_m)$ holds.
\end{remark}

The following result, required to prove the main result of this
section, holds for any ACM scheme of codimension two. Here, we
give a proof in the bihomogeneous case.

\begin{theorem} \label{nonzerodivisors}
Suppose that $Z$ is a fat point scheme in $\popo$.  If $Z$ is
ACM, then there exists $L_1,L_2 \in R$ such that $\deg L_1 =
(1,0)$ and $\deg L_2 = (0,1)$, and $L_1,L_2$ give rise to a
regular sequence in $R/\Iz$.
\end{theorem}

\begin{proof}
The Krull dimension of $R/\Iz$ is K-$\dim R/\Iz =2$. Because $Z$
is ACM, it follows that there exists a regular sequence of
length 2 in $R/\Iz$.  It is therefore sufficient to show that
the elements in the regular sequence have the appropriate
degrees.

By Lemma \ref{nonzerodivisor} there exists $L_1 \in R$ such
that $\deg L_1 = (1,0)$ and $\overline{L}_1$ is a non-zero divisor
of $R/\Iz$.  It is therefore enough to show there exists
a non-zero divisor $\overline{L}_2 \in R/(\Iz,L_1)$ with
$\deg L_2 = (0,1)$.

Let $(\Iz,L_1) = Q_1 \cap \cdots \cap Q_s$ be the primary decomposition
of $(\Iz,L_1)$ and set $\wp_i := \sqrt{Q_i}$.  We claim that
$(x_0,x_1) \subseteq \wp_i$ for each $i$.  Indeed, since $L_1$
is a non-zero divisor, we have the following exact graded
sequence:
\[
0 \longrightarrow (R/\Iz)(-1,0) \stackrel{\times L}{\longrightarrow}
R/\Iz \longrightarrow R/(\Iz,L) \longrightarrow 0.
\]
Thus, $H_{R/(\Iz,L_1)}(i,j) = H_{Z}(i,j) - H_{Z}(i-1,j)$ for all
$(i,j) \in \N^2$.  By Corollary \ref{eventualbehaviour}, if $i \gg 0$,
$H_{Z}(i,0) = H_{Z}(i-1,0)$, and hence, $H_{R/(\Iz,L_1)}(i,0) =
0$. This implies $(\Iz,L_1)_{i,0} = R_{i,0} =
\left[(x_0,x_1)^i\right]_{i,0}$. So, $(x_0,x_1)^i \subseteq Q_j$
for $i \gg 0$ and for each $j =1,\ldots,s$.  Therefore,
$(x_0,x_1) \subseteq \wp_j$ for each $j$.

The set of zero divisors of $R/(\Iz,L_1)$, denoted ${\bf Z}(R/(\Iz,L_1))$,
are precisely the elements of
\[
{\bf Z}(R/(\Iz,L_1)) = \bigcup_{i=1}^s \overline{\wp}_i.
\]
Because ${\bf k}$ is infinite, it is enough to show that
$(\wp_i)_{0,1} \subsetneq R_{0,1}$ for each $i$. If there exists
an $i \in \{1,\ldots,s\}$ such that $(\wp_i)_{0,1} = R_{0,1}$,
then $ (x_0,x_1,y_0,y_1) \subseteq \wp_i$.  But then every
homogeneous element of $R/(\Iz,L_1)$ is a zero divisor,
contradicting the fact that $Z$ is ACM. So $R/(\Iz,L_1)$ has a
non-zero divisor of degree $(0,1)$.
\end{proof}

\begin{corollary} \label{bigradedartinian}
If $Z$ is an ACM fat point scheme in $\popo$,  then the first
difference function $\Delta H_{Z}$ is the Hilbert function of a
bigraded artinian quotient of ${\bf k}[x_1,y_1]$.
\end{corollary}

\begin{proof}
Let $L_1, L_2$ be the regular sequence of Theorem \ref{nonzerodivisors}.
By making a linear change of coordinates in the $x_0,x_1$'s, and
a linear change of coordinates in the $y_0,y_1$'s, we can assume
that the $L_1 = x_0, L_2 = y_0$ give rise to a regular sequence in $R/\Iz$.

From the short exact sequences
\[
\begin{array}{cccccccccc}
0 & \rightarrow & (R/\Iz)(-1,0) &\stackrel{\times \overline{x}_0}{\rightarrow}&
&R/\Iz &\rightarrow& R/(\Iz,x_0) &\rightarrow& 0 \\
0 & \rightarrow & (R/(\Iz,x_0))(0,-1) &\stackrel{\times \overline{y}_0}
{\longrightarrow}&
&R/(\Iz,x_0) &\rightarrow& R/(\Iz,x_0,y_0) &\rightarrow& 0 \\
\end{array}
\]
it follows that $H_{R/(\Iz,x_0,y_0)}(i,j) = \Delta H_{Z}(i,j)$
for all $(i,j) \in \N^2$.  Moreover,
\[
R/(\Iz,x_0,y_0) \cong \frac{R/(x_0,y_0)}{(\Iz,x_0,y_0)/(x_0,y_0)}
\cong {\bf k}[x_1,y_1]/J
\]
where $J$ is a bihomogeneous ideal with $J \cong
(\Iz,x_0,y_0)/(x_0,y_0)$.  By using Corollary
\ref{eventualbehaviour} it follows that $\Delta H_{Z}(i,j) = 0$
if $i \gg 0$ or $j \gg 0$.  Hence ${\bf k}[x_1,y_1]/J$ is an
artinian ring.
\end{proof}

\begin{lemma} \label{rowcolumnsum}
Let $Z$ be a fat point scheme of $\popo$. Set $c_{i,j} := \Delta
H_{Z}(i,j)$. Then
\begin{enumerate}
\item[$(i)$] for every $0 \leq j \leq |\bz| -1$
\[
\alpha^*_{j+1} = \sum_{h \leq |\az| -1} c_{h,j}.
\]
where $\alpha^*_{j+1}$ is the $(j+1)$-th entry
of $\az^*$, the conjugate of the partition $\az$.
\item[$(ii)$] for every $0 \leq i \leq |\az| -1$
\[
\beta^*_{i+1} = \sum_{h \leq |\bz| -1} c_{i,h}.
\]
where $\beta^*_{i+1}$ is the $(i+1)$-th entry
of $\bz^*$, the conjugate of the partition $\bz$.
\end{enumerate}
\end{lemma}

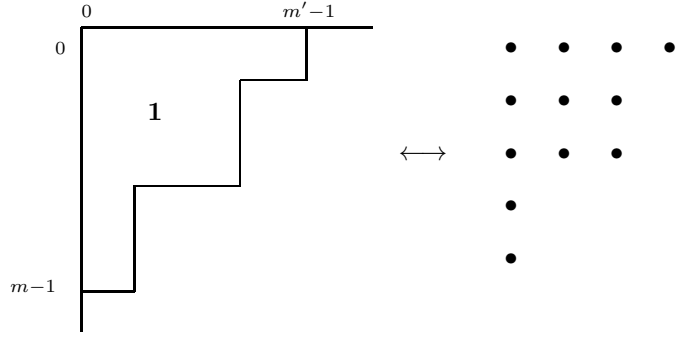
\begin{proof}
Fix an integer $j$ such that $0 \leq j \leq |\bz| -1$ and set $m
= |\az|$.  Using Theorem \ref{borderresult} and the identity
$H_{Z}(i,j) = \sum_{(h,k) \leq (i,j)} c_{h,k}$ to compute
$\alpha^*_{j+1}$ we have
\begin{eqnarray*}
\alpha^*_{j+1} & = & H_{Z}(m-1,j) - H_{Z}(m-1,j-1) \\
& = & \sum_{(h,k)\leq (m-1,j)} c_{h,k} -
\sum_{(h,k) \leq (m-1, j-1)} c_{h,k} =
\sum_{h \leq m-1} c_{h,j}.
\end{eqnarray*}
The proof for the second statement is the same.
\end{proof}

\begin{lemma}[{\cite[Lemma 6.10]{VT2}}] \label{conjugateprop}
Let $\alpha = (\alpha_1,\ldots,\alpha_n)$, $\beta = (\beta_1,\ldots,\beta_m)$,
and suppose that $\alpha, \beta \vdash s$.  If $\alpha^* = \beta$, then
\begin{enumerate}
\item[$(i)$] $\alpha_1 = |\beta|$.

\item[$(ii)$] $\beta_1 = |\alpha|$.

\item[$(iii)$]  if $\alpha' = (\alpha_2,\ldots,\alpha_n)$ and $\beta'=
(\beta_1 -1,\ldots,\beta_{\alpha_2} - 1)$, then $(\alpha')^* = \beta'$.
\end{enumerate}
\end{lemma}

\begin{theorem} \label{equivalent}
Let $Z$ be a fat point scheme in $\popo$ with Hilbert function
$H_{Z}$.  Then the following are equivalent:
\begin{enumerate}
\item[$(i)$]  $Z$ is arithmetically Cohen-Macaulay.
\item[$(ii)$] $\Delta H_{Z}$ is the Hilbert function of a bigraded
artinian quotient of ${\bf k}[x_1,y_1]$.
\item[$(iii)$] $\az^* = \bz$.
\item[$(iv)$] The set $\Ssz$ is totally ordered.
\end{enumerate}
\end{theorem}

\begin{proof}
In light of Theorem \ref{acm} and Corollary \ref{bigradedartinian}, it is
enough to prove that $(ii) \Rightarrow (iii) \Rightarrow (iv)$.

Suppose that $\Delta H_{Z}$ is the Hilbert function of a
bigraded artinian quotient of ${\bf k}[x_1,y_1]$.   Since
$\dim_{\bf k}{\bf k}[x_1,y_1]_{i,j} = 1$ for all $(i,j)$,
$\Delta H_{Z}(i,j) = 1$ or $0$.
If we write
$\Delta H_{Z}$ as an infinite matrix whose index starts from
zero, rather than one, then we have
\begin{center}
\begin{picture}(150,150)(-40,0)
\put(-60,75){$\Delta H_{Z} =$} \put(10,25){\line(0,1){115}}
\put(10,144){${\scriptstyle 0}$} \put(0,130){${\scriptstyle 0}$}
\put(86,144){${\scriptstyle m'-1}$} \put(-17,40){${\scriptstyle
m-1}$} \put(10,140){\line(1,0){110}} \put(10,40){\line(1,0){20}}
\put(30,40){\line(0,1){40}} \put(30,80){\line(1,0){40}}
\put(70,80){\line(0,1){40}} \put(70,120){\line(1,0){25}}
\put(95,120){\line(0,1){20}} \put(35,105){{\bf 1}}
\put(80,50){{\bf 0}}
\end{picture}
\end{center}
where $m = |\az|$ and $m' = |\bz|$.  By Lemma \ref{rowcolumnsum}
the number of $1$'s in the $(i-1)^{th}$ row of $\Delta H_{Z}$
for each integer $1 \leq i \leq m$ is simply the $i^{th}$
coordinate of $\bz^*$.  Similarly, the number of $1$'s in the
$(j-1)^{th}$ column of $\Delta H_{Z}$ for each integer $1 \leq j
\leq m'$ is the $j^{th}$ coordinate of $\az^*$.  Now $\Delta
H_{Z}$ can be identified with the Ferrer's diagram (see
Definition \ref{ferrers}) by associating each $1$ in $\Delta
H_{Z}$ with a dot in the Ferrer's diagram in a natural way:

\begin{picture}(150,150)(-75,10)

\put(10,25){\line(0,1){115}}
\put(10,144){${\scriptstyle 0}$}
\put(0,130){${\scriptstyle 0}$}
\put(86,144){${\scriptstyle m'-1}$}
\put(-17,40){${\scriptstyle m-1}$}
\put(10,140){\line(1,0){110}}
\put(10,40){\line(1,0){20}}
\put(30,40){\line(0,1){40}}
\put(30,80){\line(1,0){40}}
\put(70,80){\line(0,1){40}}
\put(70,120){\line(1,0){25}}
\put(95,120){\line(0,1){20}}
\put(35,105){{\bf 1}}
\put(170,50){$\bullet$}
\put(170,70){$\bullet$}
\put(170,90){$\bullet$}
\put(170,130){$\bullet$}
\put(170,110){$\bullet$}
\put(190,90){$\bullet$}
\put(190,110){$\bullet$}
\put(190,130){$\bullet$}
\put(210,90){$\bullet$}
\put(210,110){$\bullet$}
\put(210,130){$\bullet$}
\put(230,130){$\bullet$}
\put(130,90){$\longleftrightarrow$}
\end{picture}

\noindent
By using the Ferrer's diagram and Lemma \ref{rowcolumnsum} we
can calculate that $\bz = (\bz^*)^* = \az^*,$ and so $(iii)$ holds.

Now suppose that $Z$ is a fat point scheme $Z =
\{(P_{ij};m_{ij}) ~|~ 1 \leq i \leq r, 1\leq j \leq t\}$ where
$m_{ij}$ are non-negative numbers and $\az^* = \bz$.  We will
work by induction on $\beta_1 = \max\{\sum_{i=1}^r
m_{ij}\}_{j=1}^t.$

If $\beta_1 = 1$, then $Z$ is a set of $s$ distinct simple
points with $\az = (s)$ and $\bz =
(\underbrace{1,\ldots,1}_s)$.  So $Z = \{P\times Q_1,\ldots,
P\times Q_s\}$, in which case it can be easily checked that
$\Ssz = \{(1,\ldots,1)\}$, and that the set is trivially ordered.

Let us suppose that $\beta_1 > 1$ and the theorem holds for
all fat point schemes $\Y$ with $\alpha_{\Y}^* = \beta_{\Y}$,
and the first coordinate of $\beta_{\Y}$ is less than $\beta_1$.

Let $k$ be the index in $\{1,\ldots,r\}$ such that $\alpha_1
= \sum_{j=1}^t m_{kj}$.

\noindent
{\it Claim. } $m_{kj} > 0 $ for $j=1,\ldots,t$.

\noindent {\it Proof of the Claim. }  Set $l'_j =
\max\{m_{1j},\ldots,m_{rj}\}$ for $j=1,\ldots,t$.  Then $|\bz| =
l'_1 + \cdots + l'_t$.  Since $\az^* = \bz$, by Lemma
\ref{conjugateprop} $\alpha_1 = l'_1 + \cdots + l'_t$.  Now
suppose that $m_{kc} = 0$ for some $c \in \{1,\ldots,t\}$. Since
$l'_j \geq m_{kj}$ for each $j=1,\ldots,r$, we would then have
\begin{eqnarray*}
\alpha_1  =  l'_1 + \cdots + l'_t
 & > & l'_1 + \cdots + \hat{l}'_c + \cdots + \l'_t \\
 & \geq & m_{k1} + \cdots + \hat{m}_{kc} + \cdots + m_{kt} \\
 & = & m_{k1} + \cdots + m_{kc} + \cdots m_{kt} = \alpha_1
\end{eqnarray*}
where \quad$\hat{ }$\quad means the number is omitted. Because of this
contradiction, the claim holds. \hfill $\diamond$

Let $\Y = \{(P_{ij};m'_{ij}) ~|~ 1 \leq i \leq r, 1 \leq j \leq
t \}$ be the subscheme of $Z$ where
\[
m'_{ij} = \left\{
\begin{array}{ll}
m_{ij} & i \neq k \\
m_{kj}-1 & i = k
\end{array}
\right.
\]
By the claim $m_{kj} - 1 \geq 0$ for all $j = 1,\ldots,t$.
Let $\beta$ be the first coordinate of $\beta_{\Y}$.  Then $\beta < \beta_1$.
In fact, for each $j = 1,\ldots,t$, we have
\[
\sum_{i=1}^r m'_{ij} = m'_{kj} + \sum_{i \neq k} m_{ij} =
\left(\sum_{i=1}^r m_{ij} \right) - 1.
\]
Furthermore, if $\az = (\alpha_1,\ldots
\alpha_m)$ and $\bz = (\beta_1,\ldots,\beta_{m'})$, then
from our construction $\Y$ we have $\alpha_{\Y} = (\alpha_2,\ldots,\alpha_m)$
and $\beta_{\Y} = (\beta_1 - 1,\ldots,\beta_{\alpha_2} -1)$.
By Lemma \ref{conjugateprop}, $\alpha_{\Y}^* = \beta_{\Y}$, and so by
induction $\mathcal{S}_{\Y}$ is totally ordered.

The set $\Ssz$ is now obtained from $\mathcal{S}_{\Y}$ by adding
the tuple $(m_{k1},\ldots,m_{kt})$.  Moreover, this element is
larger than every other element of $\mathcal{S}_{\Y}$ with
respect to our ordering, so $\Ssz$ is totally ordered, as
desired.
\end{proof}

\begin{corollary}
If $Z$ is a scheme of fat points whose support is on a line,
then $Z$ is ACM.
\end{corollary}

\begin{proof}
It easy to check that either the set $\Ssz$ is
totally ordered, or $\az^* = \bz$.
\end{proof}

\begin{corollary} \label{fatpoints}
If $Z$ is an ACM scheme of fat points with $\az =
(\alpha_1,\ldots,\alpha_m),$ then the Hilbert function of $Z$ is
\begin{eqnarray*}
H_{Z} & = &  \bmatrix
1& 2 & \cdots  & \alpha_{1}-1 & \alpha_1 & \alpha_1 & \cdots \\
1& 2 &\cdots  & \alpha_{1}-1 & \alpha_1 & \alpha_1 & \cdots \\
\vdots & \vdots && \vdots  & \vdots &\vdots&\ddots
\endbmatrix
+ \bmatrix
0& 0 & \cdots  & 0 & 0 & 0 & \cdots\\
1& 2 & \cdots  & \alpha_{2}-1 & \alpha_2 & \alpha_2 & \cdots \\
1& 2 & \cdots  & \alpha_{2}-1 & \alpha_2 & \alpha_2 & \cdots \\
\vdots & \vdots && \vdots  & \vdots &\vdots&\ddots
\endbmatrix \\
& & + \cdots + \bmatrix
0& 0 & \cdots  & 0 & 0 & 0 &\cdots\\
\vdots& \vdots& & \vdots & \vdots & \vdots \\
0& 0 & \cdots  & 0 & 0 & 0 &\cdots\\
1& 2 & \cdots  & \alpha_{m}-1 & \alpha_m & \alpha_m  & \cdots \\
1& 2 & \cdots  & \alpha_m-1 & \alpha_m & \alpha_{m} & \cdots \\
\vdots & \vdots & &\vdots  & \vdots &\vdots&\ddots
\endbmatrix.
\end{eqnarray*}
\begin{proof}
Use Theorem \ref{acm} and Remark \ref{az=sz}.
\end{proof}
\end{corollary}

From the above corollary, we see that if the fat point scheme
$Z$ in $\popo$ is ACM, then the entire Hilbert function of $Z$
can be determined from the tuple $\az$. This contrasts with the
main result of the previous section where we showed that for a
general fat point scheme in $\popo$, most, but not all,  of the
values of the Hilbert function can be determined from the tuples
$\az$ and $\bz$.

In fact, if $Z$ is an ACM fat point scheme in $\popo$, we can
even compute the Betti numbers in the minimal free resolution of
$\Iz$ directly from the tuple $\az$.  To state our result, we first
develop some suitable notation.

Let $Z$ be an ACM scheme of fat points and let
$\az=(\alpha_1,\ldots,\alpha_m)$ be the tuple associated to $Z$.
Define the following two sets from $\az$:
\begin{eqnarray*}
C_{Z} & := & \left\{(m,0),(0,\alpha_1)\right\} \cup
\left\{(i-1,\alpha_i) ~|~ \alpha_i - \alpha_{i-1} < 0\right\} \\
V_{Z} & := & \left\{ (m,\alpha_m) \right\} \cup \left\{
(i-1,\alpha_{i-1}) ~|~ \alpha_i-\alpha_{i-1} < 0 \right\}.
\end{eqnarray*}
We take $\alpha_{-1} = 0$.  With this notation, we have

\begin{theorem} \label{bettinumbers}
Suppose that $Z$ is an ACM set of fat points in $\popo$ with
$\az = (\alpha_1,\ldots,\alpha_m)$.  Let $C_{Z}$ and $V_{Z}$ be
constructed from $\az$ as above.  Then the bigraded minimal free
resolution of $\Iz$  is given by
\[
0 \longrightarrow \bigoplus_{(v_1,v_2) \in V_{Z}} R(-v_1,-v_2)
\longrightarrow \bigoplus_{(c_1,c_2) \in C_{Z}} R(-c_1,-c_2)
\longrightarrow \Iz \longrightarrow 0.\]
\end{theorem}

\begin{proof}
Using Theorem \ref{acm}, it can be verified that the tuples in
the set $C_{Z}$ are what \cite{GuMaRa1} defined to be the {\it
corners} of $\Delta H_{Z}$, and the elements in $V_{Z}$ are
precisely the {\it vertices} of $\Delta H_{Z}$.  The conclusion
now follows from Theorem $4.1$ in \cite{GuMaRa1} .
\end{proof}


\section{Special configurations of ACM fat points}

Theorem \ref{equivalent} enables us to identify the ACM fat point
schemes directly from the tuples $\az$ and $\bz$, or from the set
$\Ssz$.  In this section, we use these characterizations to
investigate ACM fat point schemes which have some extra
conditions on the multiplicities of the points.  We show that
some
 special configurations of ACM fat point schemes
can occur only if the support of the scheme has some
specific properties.

\begin{remark}\label{rednonacm}
By Theorem $2.12$ and Theorem $4.1$ in \cite{GuMaRa1}, we can
deduce that $\X$ is not an ACM scheme if and only if there exist two points
$P_{11}=[a_1:a_2]\times[b_1:b_2]$ and $P_{22}=
[c_1:c_2] \times[d_1:d_2]$ of $\X$
with $a_i,b_i,c_i,d_i \in {\bf k}$ such that
$P_{12}=[a_1:a_2]\times[d_1:d_2]$ and
$P_{21}=[c_1:c_2] \times[b_1:b_2] \not\in \X$.
\end{remark}

\begin{proposition} \label{suppci}
If $Z$ is an ACM fat point scheme, then $\supp(Z)$ is ACM.
\end{proposition}

\begin{proof}
Let us suppose that $\supp(Z)$ is not ACM.   Then by Remark
\ref{rednonacm}, in $\Ssz$ we can find tuples of type:
\[
(*,1,*,0,*),(*,0,*,1,*)
\]
\noindent that are incomparable. Therefore, by Theorem
\ref{equivalent}, $Z$ is not ACM.
\end{proof}

\begin{remark}  \label{ciremark}
Theorem 1.2 of \cite{GuMaRa1}
showed that for any saturated bihomogeneous
ideal $I \subseteq R$ of height two,
the minimal generating set for $I$ must
contain exactly one form of degree $(m,0)$ for
some $m$, and one form of degree $(0,n)$ for some $n$.
If $F \in I$ is the form of degree $(m,0)$, then
$F \in {\bf k}[x_0,x_1] \subseteq R$, and thus $F$ can
be written as the product of $(1,0)$ forms.  Similarly,
the form of degree $(0,n)$ can be written
as a product of forms of degree $(0,1)$.  Thus, following
Remark 1.3 of \cite{GuMaRa1}, we shall call a set
of points $\X$ a {\it complete intersection} if $\Ix = (F,G)$
where $\deg F = (m,0)$ and $\deg G = (0,n)$.
\end{remark}

We now describe the support of the ACM fat point schemes which are
{\it homogeneous}, i.e., all the nonzero multiplicities are equal.

\begin{theorem} \label{completeintersection}
Fix a positive integer $m \geq 2$, and let $Z$ be a homogeneous fat
point scheme of $\popo$ with all the nonzero multiplicities
equal to $m$. Then $Z$ is ACM if and only if $\supp(Z)$ is a
complete intersection.
\end{theorem}

\begin{proof}
If $\supp(Z)$ is a complete intersection, then $Z$ is ACM by
Corollary 2.5 of \cite{Gu}.

Conversely, suppose that $Z$ is ACM, and thus, $\Ssz$ is totally
ordered by Theorem \ref{equivalent}. Because $Z$ is ACM, from
Proposition \ref{suppci}, $\supp(Z)$ must also be ACM.

Suppose that $\supp(Z)$ is not a complete intersection. This
implies that $Z$ contains a subscheme of type
\begin{eqnarray*}
\Y &=  &\{(P_{i_{1}j};m_{i_{1}j}) ~|~  m_{i_{1}j}=m ~\text{for}~ j=1,\ldots,t\}
\cup \\
& &
\left\{(P_{i_{2}j};m_{i_{2}j})
~\left|~
\begin{array}{ll}
m_{i_2j} = m &  j=1,\ldots,h \mbox{ with $h<t$} \\
m_{i_2j} = 0 &\mbox{otherwise}.
\end{array}
\right.
\right\}
\end{eqnarray*}
But then in $\Ssz$ we can find three tuples of
the form
\[
\{(\underbrace{m,\ldots,m}_t),({\underbrace{m,\ldots,m,}_h}\underbrace{0,\ldots,0}_{t-h}),(\underbrace{m-1,\ldots,m-1}_t)\}.
\]
But then $\Ssz$ is not totally ordered, which is a contradiction.
\end{proof}

\begin{remark}
Homogeneous schemes with all $m_{ij} = 2$ have been further
investigated by the first author in \cite{Gu}.
\end{remark}

\begin{definition}
A fat point scheme $Z$ in $\popo$ is called an {\it almost
homogeneous fat point scheme} if all the non-zero multiplicities
of $Z$ are either $m$ or $m-1$ for some integer $m >0$.
\end{definition}

We now recall a definition first given in \cite{Gu}.

\begin{definition}
Let $Z = \{(P_{ij};m_{ij}) ~|~ 1\leq i \leq r, 1\leq j \leq t\}$
be a fat point scheme.  The scheme $Z$ is called a {\it
quasi-homogeneous scheme of fat points} if there exist $r$
integers $t = t_1 \geq t_2 \geq \cdots \geq t_r\geq 1$ such that
\[
m_{ij} = \left\{
\begin{array}{ll}
m & j=1,\ldots,t_i \\
m-1 & j = t_{i+1},\ldots,t_1
\end{array}
\right.
.\]
\end{definition}

\begin{remark}
Note that if $Z$ is a quasi-homogeneous scheme and $m \geq 2$,
then $\supp(Z)$ is the complete intersection $\{P_{ij} ~|~ 1
\leq i \leq r, ~1 \leq j \leq t\}$.   If $m=1$, then a
quasi-homogeneous scheme of fat points is an ACM scheme of
simple points.  However, if $m=1$, then the support is not a
complete intersection.
We also observe that any quasi-homogeneous fat point
scheme is also an almost homogeneous fat point scheme for any $m$.
\end{remark}

\begin{remark}
If $Z$ is a quasi-homogeneous fat point scheme, then $Z$ is ACM by
Corollary $2.6$ in \cite{Gu} .
\end{remark}

Since $\popo$ is isomorphic to the quadric surface $\Q \subseteq \pr^3$,
using Remark \ref{ciremark}, we can draw fat point schemes
on $\Q$ as subschemes whose support is contained in the
intersection of lines of the two rulings of $\Q$.
For example,
if $P_{ij} = R_i \times Q_j \in \popo$, then the fat point scheme
$\Z = \{(P_{11};4),(P_{12};2),(P_{22};3)\}$
can be visualized as

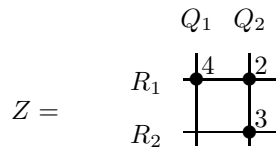
\begin{picture}(100,60)(-70,5)
\put(10,20){$Z = $} \put(80,10){\line(0,1){35}}
\put(100,10){\line(0,1){35}}
\put(74,55){$Q_1$} \put(94,55){$Q_2$}
\put(75,15){\line(1,0){35}}
\put(75,35){\line(1,0){35}}
\put(55,11){$R_2$}
\put(55,31){$R_1$}
\put(80,35){\circle*{5}}\put(82,37){4}
\put(100,35){\circle*{5}} \put(102,37){2}
\put(100,15){\circle*{5}} \put(102,17){3}
\end{picture}

\noindent
where a dot represents a point in the support, and the number
its multiplicity.

\begin{theorem} \label{quasihomogeneous}
Let $Z$  be a fat point scheme. If $Z$ is an ACM almost
homogeneous fat point scheme with $m \geq 4$, then $Z$ is a
quasi-homogeneous scheme of fat points. In particular, the
support of $Z$ is a complete intersection.
\begin{proof} Suppose that $Z$ is an ACM  almost homogeneous fat point scheme.

\noindent {\it Claim. } $\supp(Z)$ is a complete intersection.

\noindent {\it Proof of the Claim. }  For a contradiction,
suppose $\supp(Z)$ is not a complete intersection.  Since
$\supp(Z)$ is contained within a complete intersection, we can
find a point $P_{ij} = R_i \times Q_j \not\in \supp(Z)$ but
$P_{i'j} = R_{i'} \times Q_j$ and $P_{ij'} = R_i \times Q_{j'}$
in $\supp(Z)$. So $Z$ contains the following subscheme

\begin{centering}
\begin{picture}(100,60)(-100,45)
\put(80,50){\line(0,1){30}} \put(100,50){\line(0,1){30}}
\put(74,95){$Q_{j'}$} \put(94,95){$Q_j$}
\put(75,55){\line(1,0){30}} \put(75,75){\line(1,0){30}}
\put(55,51){$R_{i}$} \put(55,71){$R_{i'}$}
\put(80,55){\circle*{5}} \put(82,57){$b$}
\put(80,75){\circle*{5}} \put(82,77){$c$}
\put(100,75){\circle*{5}} \put(102,77){$a$}
\put(100,55){\circle*{5}} \put(102,57){$0$}
\end{picture}\end{centering}

\noindent
where $a,b,$ and $c$ denote the multiplicities of $R_{i'} \times Q_j,
R_i \times Q_{j'}$ and $R_{i'} \times Q_{j'}$ respectively, and $0$
denotes the absence of the point $R_i \times Q_j$.

We observe that the tuples $(*,c,*,a,*)$ and $(*,b,*,0,*)$ are in
$\Ssz$ with $c$ and $b$ in the $j^{'th}$ spot and the $a$ and $0$
in the $j^{th}$ spot, and where $*$ denotes the other unknown
numbers in the tuple.  Because $Z$ is ACM, $\Ssz$ is totally
ordered, so $m \geq c \geq b \geq m-1$.

We see that $c$ can be either $c > b$ or $c = b$.  If $c > b$,
then $c = m$ and $b = m -1$.  But then the tuple $(*,
m-2,*,a-2,*)$ is also in $\Ssz$ with $a -2 \geq (m-1)-2 > 0$
because $m \geq 4$.  But then $\Ssz$ is not totally ordered
because the tuples $(*,b,*,0,*)$ and $(*,c-2,*,a-2,*)$ are
incomparable.

Similarly, if $c = b$, then the tuple $(*,c-1,*,a-1,*)$ is in
$\Ssz$ with $b > c-1$, but $a-1 > 0$, contradicting the fact that
$\Ssz$ is totally ordered.
So, the support of $Z$ must be a complete
intersection.\hfill$\diamond$

Because of the claim, we can consider subschemes of $Z$ that
consist of the following four points: $P_{ij} = R_i \times Q_j$,
$P_{i'j} = R_{i'} \times Q_j$, $P_{ij'} = R_i \times Q_{j'}$,
and $P_{i'j'} = R_{i'} \times Q_{j'}$.  Now no such subscheme
will have the form

\begin{picture}(100,60)(-100,45)

\put(80,50){\line(0,1){30}}
\put(120,50){\line(0,1){30}}

\put(74,95){$Q_{j'}$}
\put(114,95){$Q_j$}

\put(75,55){\line(1,0){50}}
\put(75,75){\line(1,0){50}}

\put(55,51){$R_{i}$}
\put(55,71){$R_{i'}$}

\put(80,55){\circle*{5}}
\put(82,57){$m$}

\put(80,75){\circle*{5}}
\put(82,77){$m-1$}

\put(120,75){\circle*{5}}
\put(122,77){$m$}

\put(120,55){\circle*{5}}
\put(122,57){$m-1$}
\end{picture}

\noindent because such a subscheme would contradict the fact that
$\Ssz$ is totally ordered. So, if we write only the multiplicities
of the points, then the scheme $Z$ must have the form
\[
\begin{matrix}
m & m & \cdots & m & m & m\\
\vdots &\vdots &  &\vdots & \vdots&\vdots\\
m & m & \cdots & m & m & m \\
m & m & \cdots & m & m & m-1\\
m & m & \cdots & m & m-1 & m-1 \\
\vdots &\vdots &  & \vdots&\vdots& \vdots\\
m & m-1 & \cdots &m-1&  m-1 & m-1 &\end{matrix}
\]
that is, $Z$ is a quasi-homogeneous scheme of fat points.
\end{proof}
\end{theorem}
\begin{example}
One can check that the following scheme

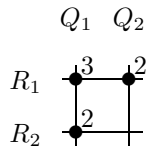
\begin{picture}(100,60)(-100,45)

\put(80,50){\line(0,1){30}}
\put(100,50){\line(0,1){30}}

\put(74,95){$Q_{1}$}
\put(94,95){$Q_2$}

\put(75,55){\line(1,0){30}}
\put(75,75){\line(1,0){30}}

\put(55,51){$R_{2}$}
\put(55,71){$R_{1}$}

\put(80,55){\circle*{5}}
\put(82,57){$2$}

\put(80,75){\circle*{5}}
\put(82,77){$3$}

\put(100,75){\circle*{5}}
\put(102,77){$2$}

\end{picture}

\noindent is an almost homogeneous fat point scheme that is also
ACM. However, the support is not a complete intersection.  So the
hypothesis $m \geq 4$ is needed in the above theorem.
\end{example}


\end{document}